\documentclass{crmp-l}
\copyrightinfo{2006}{American Mathematical Society}

\newtheorem{theorem}[subsection]{Theorem}
\newtheorem{proposition}[subsection]{Proposition}
\newtheorem{lemma}[subsection]{Lemma}
\newtheorem{corollary}[subsection]{Corollary}

\def\zp{\mathbb{Z}_p}
\def\ve{\varepsilon}

\begin{document}
\title{Additive properties of product sets in fields of prime
order}

\author{Glibichuk A.~A.}

\author{Konyagin S.~V.}
\thanks{This research was carried out while the authors were
visitors of le Centre de Recherches Math\'ematiques. It is our
pleasure to thank CRM for its hospitality and the Clay Institute for
a generous support of our visit.}

\maketitle

\section{Introduction} \label{sec1}

Let $p>2$ be a prime, $\zp$ be the field of the residues modulo $p$, and
$\zp^*$ be the multiplicative group of $\zp$. Thus,
$\zp^*=\zp\setminus\{0\}$.
For sets
$X\subset\zp$, $Y\subset\zp$, and for
(possibly, partial) binary operation $*:\zp\times\zp\to\zp$ we denote
$$X*Y=\{x*y:\,x\in X,y\in Y\}.$$
Usually we will write $XY$ instead of $X*Y$ if $*$ is the multiplication.
Also, for an element $\xi\in\zp$ denote
$$\lambda\ast A=\{\lambda\}A.$$
For a set $X\subset\zp$ and $k\in\mathbb{N}$ we denote
$$kX=\{x_1+\dots+x_k:\,x_1,\dots,x_k\in A\},$$
$$X^k=\{x_1\dots x_k:\,x_1,\dots,x_k\in A\}.$$
A set $X$ is called a basis (an additive basis) of order $k$ if $kA=X$.
Observe that any basis of order $k$ is also a basis of any order $k'>k$.
A general problem to be discussed in this paper is the following:
given $p, t\leq p, n, N$. Is it true that for any set $A\subset\zp$
of cardinality $\geq t$, the set $A^n$ is the basis of order $N$?

The situation is clear for $n=1$.
Due to Cauchy--Davenport theorem (\cite{TV}, Theorem 5.4) for any sets
$X_1,\dots,X_N\subset\zp$ we have
$$|X_1+\dots+X_N|\geq\min(|X_1|+\dots+|X_N|-N+1,p).$$
Therefore, any set $A$ with $|A|-1\geq(p-1)/N$ is a basis of order $N$.
On the other hand, if $t\in\mathbb{N}$ and $t<1+(p-1)/N$, it is easy
to see that a set $A=\{0,\dots,t-1\}$ satisfy the conditions
$|A|=t$ and $NA\neq\zp$.

In the case $n=2$ some useful information can be obtained
by using exponential sums. Some facts related to the harmonic analyisis
in $\zp$ can be found in (\cite{TV}, Chapter 4).

It is known that for fixed $k\in\mathbb{N}$ and $\ve>0$ a random subset
of $\zp$ of cardinality $>p^{\ve+1/k}$ is a basis of order $k$ with
a large probability (tending to $1$ as $p\to\infty$). Thus, if
a large set is not a basis of small order, it has a special additive structure.
We can believe (and have many confirmations) that this cannot hold for sets
possessing a special multiplicative structure. If a set does,
in many cases nontrivial estimates for exponential sums over the set and
additive properties of the set are known. Probably, the simplest result
of this type is the following (\cite{Vin}, Chapter VI, Problem 8,$\alpha$).
\begin{proposition}\label{prop1} If $X,Y\subset\zp$, $a\in\zp^*$, then
$$\left|\sum_{x\in X, y\in Y}\exp(2\pi iaxy/p)\right|\leq\sqrt{p|X||Y|}.$$
\end{proposition}

Using Proposition\ref{prop1} and standard technique, one can easily deduce
that for any $\ve>0$ there is $N=N(\ve)$ such that for any $A\subset\zp$
with $|A|>p^{1/2+\ve}$ we have
\begin{equation}\label{2,N}
NA^2=\zp.
\end{equation}
In particular, it is known that (\ref{2,N}) holds for
$|A|>p^{3/4}+1$ and $N=3$.

However, we do not see a way to prove (\ref{2,N}) with a bounded $N$
via exponential sums since there is no estimate for exponential sums
under a weaker restriciton
\begin{equation}\label{sqrt{p}}
|A|>\sqrt p
\end{equation}
essentially better than the trivial estimate $|A|^2$. In \cite{Gli},
by using combinatorial arguments, it has been proven that $8XY=\zp$
for any $X\subset\zp$, $Y\subset\zp$ provided that $|X||Y|>p$
and either $Y=-Y$, or $Y\cap(-Y)=\emptyset$. This easily implies
(\ref{2,N}) in the case (\ref{sqrt{p}}) with $N=16$ (see Section 2).
Restrictiontion (\ref{sqrt{p}}) is essentially sharp.
Clearly, if $f(p)=o(p^{1/2})$ as $p\to\infty$, then the condition
$|A|\leq f(p)$ cannot guarantee that $A^2$ is a basis of a fixed order
as can be seen by taking $A=\{1,\dots,[f(p)]\}$.

By the same reasons, for any fixed $n\in\mathbb{N}$ and $f(p)=o(p^{1/n})$
as $p\to\infty$, the condition $|A|\leq f(p)$ cannot guarantee that $A^n$
is a basis of a fixed order.

The estimates of exponential sums established in \cite{BGK} (Theorem 5)
clearly imply that for any $\delta>0$ there are $n(\delta)$ and $N(\delta)$
such that
\begin{equation}\label{n,N}
NA^n=\zp
\end{equation}
for any set $A\subset\zp$ with $|A|>p^{\delta}$; moreover, for $\delta<1/2$
we have
$$n\leq\delta^{-C},\quad N\leq\exp\left(\delta^{-C}\right)$$
where $C$ is a constant not evaluated in \cite{BGK}. It was naturally to ask
whether one could get sharper bounds for $n$ and $N$. The main result of this
paper is the following.
\begin{theorem}\label{th1} There exists a constant $C$ such that for
any integer $n>1$, any numbers $\ve\in(0,n)$, $\delta\geq 1/(n-\ve)$,
any prime $p$, and any set $A\subset\zp$ with $|A|>p^\delta$ we have
(\ref{n,N}) with
\begin{equation}\label{estN}
N\leq C4^n\log(2+1/\ve).
\end{equation}
\end{theorem}
The restriction for $n$ is essentially best possible: we have seen that
in general one could not take $n<1/\delta$. Also, $N$ should grow
at least as an exponential function of $1/\delta$ as $\delta\to0$.
This can be demonstrated by a simple example $A=\{0,1\}$. Then
$|A|>p^\delta$ if $p<2^{1/\delta}$. Next, $A^n=A$ for all $n$, and
(\ref{n,N}) holds only for $N\geq p-1$.
\begin{corollary}\label{cor1} There exists a constant $C$ such that for
any $\delta>0$, any prime $p$, and any subgroup $A$ of $\zp^*$ of cardinalty
$>p^\delta$ we have $NG=\zp$ with $N\leq C4^{1/\delta}$.
\end{corollary}
To prove Corollary \ref{cor1} it is enough to take a large $n$ and to observe
that $G^n=G$.

In the case $n>2$ we can not prove that (\ref{n,N}) holds for any $A$ with
$|A|>p^{1/n}$ and for some $N=N(n)$.

We will get some preliminary results in Sections 2--5 and prove Theorem \ref{th1}
in Section 6.

\section{On additive properties of a product of two sets} \label{sec2}

\begin{lemma}\label{lem1} If $A\subset\zp$,  $B\subset\zp$, and
$|A|\cdot\lceil|B|/2\rceil>p$ then $8AB=\zp$.
\end{lemma}
{\bf The proof of Lemma \ref{lem1}.} We split $B$ into symmetric and
antisymmetric parts:
$$B_1=\{b\in B:\,-b\in B\},\quad B_2=\{b\in B:\,-b\not\in B\}.$$
Then $|B_i|\geq \lceil|B|/2\rceil$ for $i=1$ or $i=2$. By \cite{Gli},
$8AB_i=\zp$. Hence, $8AB=\zp$, as required.
\begin{lemma}\label{lem2} If $A\subset\zp$,  $B\subset\zp$, and
$|A||B|> p$ then $16AB=\zp$.
\end{lemma}
{\bf The proof of Lemma \ref{lem2}.} By Cauchy--Davenport theorem,
we have
$$|2B|\geq\min(2|B|-1,p).$$
If $2B=\zp$, there is nothing to prove. Indeed, take $a\in A\cap\zp^*$.
Then $a*(B+B)=\zp$. Therefore, $AB+AB=\zp$ and $16AB=\zp$.
Consider the case $|2B|\geq2|B|-1$. Then $\lceil|2B|/2\rceil\geq|B|$.
As we have shown, $8A(2B)=\zp$; therefore, $16AB=\zp$.

\section{Main lemmata} \label{sec3}

For sets $X,Y\subset\zp$, $|Y|>1$, we denote
$$Q[X,Y]=\frac{X-X}{Y-Y}.$$
We will use the following observation.
\begin{lemma}\label{lem3} Let $\xi\in\zp$. Then $\xi\in Q[X,Y]$
if and only if $|X+\xi* Y|<|X||Y|$.
\end{lemma}
{\bf The proof of Lemma \ref{lem3}.} Consider the mapping $F:X\times
Y$ onto $X+\xi* Y$ defined as $F(x,y)=x+\xi y$. $F$ is not an
injection if and only if $|X+\xi* Y|<|X||Y|$. On the other hand, the
condition that $F$ is not an injection means that there are
$x_1,x_2\in X$, $y_1,y_2\in Y$ such that $F(x_1,y_1)=F(x_2,y_2)$, or
$\xi=(x_1-x_2)/(y_2-y_1)\in(X-X)/(Y-Y)$. This completes the proof of
the lemma.

For $X=Y$ Lemma \ref{lem3} is Lemma 2.50 from \cite{TV}.

The key ingredient for the proof of Theorem 2 is the following lemma
based on a technology developed by T.~Tao and V.~Vu (\cite{TV}, section 2.8).
\begin{lemma}\label{lem4} If $X,Y\subset\zp$, $a\in\zp^*$, $|Y|>1$, and
$Q[X,Y]\neq\zp$ then
$$|2XY-2XY+a*Y^2-a*Y^2|\geq|X||Y|.$$
\end{lemma}
{\bf The proof of Lemma \ref{lem4}.} By the conditions of the lemma,
$Q[X,Y]\neq\zp$ and $Q[X,Y]\neq\emptyset$. Thus, there exists
$\xi\in Q[X,Y]$ such that
$\xi+a\not\in Q[X,Y]$. The first condition implies
$\xi=(x_1-x_2)/(y_1-y_2)$ for some $x_1,x_2\in X$, $y_1,y_2\in Y$,
and from the second condition, applying Lemma 1, we get
$$|\{x+((x_1-x_2)/(y_1-y_2)+a)y:
x\in X,y\in Y\}|=|X|||Y|.$$
Multiplying by $y_1-y_2$ we get
$$|\{x(y_1-y_2)+(x_1-x_2)y+(y_1-y_2)ay:
x\in X,y\in Y\}|=|X|||Y|.$$
But
$$x(y_1-y_2)+(x_1-x_2)y+(y_1-y_2)ay
\in 2XY-2XY+a*Y^2-a*Y^2,$$
and the proof is complete.

\begin{lemma}\label{lem5} (\cite{BKT}, \cite{BK}). Let $X\subset\zp$,
$Y\subset\zp$, $G\subset\zp^*$, and $G\neq\emptyset$. Then there exists
$\xi\in G$ such that
$$|X+\xi* Y|\geq\frac{|X||Y||G|}{|X||Y|+|G|}.$$
\end{lemma}
\begin{lemma}\label{lem6} If $\xi\in Q[X,Y]$, then
$$|2XY-2XY|\geq|X+\xi*Y|.$$
\end{lemma}
{\bf The proof of Lemma \ref{lem6}.} By the condition on $\xi$, there are
elements $x_1,x_2\in X$ and $y_1,y_2\in Y$ such that $y_1\neq y_2$ and
\begin{equation}\label{equaxi}
x_1-x_2=(y_1-y_2)\xi.
\end{equation}
Denote
$$S=(y_1-y_2)*(X+\xi*Y).$$
Clearly,
$$|S|=|X+\xi* Y|.$$
Moreover, any element $s\in S$ can be written as
$$s=(y_1-y_2)x+(y_1-y_2)\xi y,\quad x\in X, y\in Y.$$
Plugging in (\ref{equaxi}) into the last equality, we get
$$s=x(y_1-y_2)+(x_1-x_2)y.$$
Thus, $S\subset 2XY-2XY$, and the proof of the lemma is complete.

\begin{corollary}\label{cor3} If $X,Y\subset\zp$, $a\subset\zp^*$,
and $|Y|>1$, then
$$|2XY-2XY+a*Y^2-a*Y^2|
\geq\frac{|X||Y|(p-1)}{|X||Y|+p-1}.$$
\end{corollary}
{\bf The proof of Lemma \ref{cor3}.} It suffices to consider the case

$Q[X,Y]=\zp$. By Lemma \ref{lem5}, there is $\xi\in\zp^*$ such that
$$|X+\xi* Y|\geq \frac{|X||Y|(p-1)}{|X||Y|+p-1}.$$
Since $\xi\in Q[X,Y]$, we can use Lemma \ref{lem6}:
$$|2XY-2XY+a*Y^2-a*Y^2|
\geq|2XY-2XY|\geq |X+\xi* Y|,$$
and we are done.
\begin{corollary}\label{cor4} If $|Y|>1$ then
$$|3Y^2-3Y^2|\geq\frac{|Y|^2(p-1)}{|Y|^2+p-1}.$$
\end{corollary}
For the proof it suffices to take $a=1$ and to use Corollary \ref{cor3}.
\begin{corollary}\label{cor5} If $|Y|>1$, $X=KY^k-KY^k$, then
$$|(4K+1)Y^{k+1}-(4K+1)Y^{k+1}|
\geq\frac{|X||Y|(p-1)}{|X||Y|+p-1}.$$
\end{corollary}
For the proof it suffices to take $y\in Y$, $y\neq0$, $a=y^{k-1}$
and to use Corollary \ref{cor3}.

\section{Ruzsa triangle inequality and its corollaries} \label{sec4}

The following nice result belongs to I.~Ruzsa (see \cite{TV}, Lemma 2.6).
\begin{lemma}\label{lem7} For any subsets $X,Y,Z$ of $\zp$ we have
$$|X||Y-Z|\leq|X-Y||X-Z|.$$
\end{lemma}
Lemma is called Ruzsa triangle inequality since it can be reformulated as
follows: the binary function $\rho$ defined for nonzero subsets of $\zp$ as
$$\rho(X,Y)=\log(|X-Y|^2/(|X||Y|))$$
satisfies the triangle inequality. The lemma is stated for subsets of an
arbitrary abelian group, but we need it only for subsets of $\zp$.

\begin{corollary}\label{cor6} If $X,Y\subset\zp$, then
$$|X+Y|\geq|X|^{1/2}|Y-Y|^{1/2}.$$
\end{corollary}
{\bf The proof of Corollary \ref{cor6}.} By Lemma \ref{lem7}, we have
$$|X||(-Y)-(-Y)|\leq|X+Y||X+Y|,$$
and we are done.
\begin{corollary}\label{cor7} If $X\subset\zp$, $k\in\mathbb{N}$, then
\begin{equation}\label{applRuz}
|kX|\geq|X|^{2^{1-k}}|X-X|^{1-2^{1-k}}.
\end{equation}
\end{corollary}
{\bf The proof of Corollary \ref{cor7}.} We use induction on $k$.
For $k=1$ the result is trivial. Now we assume that $k>1$, (\ref{applRuz})
is true for $k-1$, and we will prove it for $k$. By Corollary \ref{cor6},
$$|kX|=|(k-1)X+X|
\geq|(k-1)X|^{1/2}|X-X|^{1/2},$$
and using the induction supposition completes the proof of the corollary.

\section{Some inequalities} \label{sec5}

In this section we will establish some lower estimates for
$|KA^k|$ and $|KA^k-KA^k|$ where $A\subset\zp$ with
\begin{equation}\label{restr}
|A|\geq5.
\end{equation}
We construct the sequence of sets: $A_1=A$,
$$A_k=N_kA^k-N_kA^k,\,N_k=\frac5{24}4^k-\frac13\,(k\geq2).$$
So,
$$A_2=3A^2-3A^2,\,A_3=13A^3-13A^3,\dots$$
Corollaries \ref{cor4} and \ref{cor5} show that for $k\geq2$
\begin{equation}\label{inden}
|A_k|\geq\frac{|A||A_{k-1}|(p-1)}{|A||A_{k-1}|+p-1}.
\end{equation}
We can deduce from (\ref{inden}) an explicit lower bound for $|A_k|$.

\begin{lemma}\label{lem8} For any $k$ and $0\leq U\leq |A|^k$ we have
$$|A_k|\geq U-\frac54\frac{U^2}{p-1}.$$
\end{lemma}
{\bf The proof of Lemma \ref{lem8}.} Observe, that for any $u>0$ we have
$$\frac{u}{1+u/(p-1)}\geq u(1-u/(p-1));$$
the inequality can be rewritten as
\begin{equation}\label{basic}
\frac{u(p-1)}{u+p-1}\geq u-u^2/(p-1).
\end{equation}
We use induction on $k$. For $k=1$ the assertion of Lemma \ref{lem8} is obvious.
Assume that it holds for $k-1\geq1$ and prove it for $k$. By the induction supposition,
$$|A_{k-1}|\geq V:=\frac U{|A|}-\frac54\frac{U^2}{|A|^2(p-1)}.$$
If $V<0$ then also
$$U-\frac54\frac{U^2}{p-1}<0,$$
and the assertion of the lemma is trivial. If $V\geq0$, then,
applying (\ref{inden}) and (\ref{basic}), we have
\begin{eqnarray*}
|A_k|\geq\frac{|A||A_{k-1}|(p-1)}{|A||A_{k-1}|+p-1}
\geq\frac{|A|V(p-1)}{|A|V+p-1}\\
\geq U-\frac54\frac{U^2}{|A|(p-1)}-\frac{U^2}{p-1}=
U-\frac{U^2}{p-1}\left(1+\frac5{4|A|}\right),
\end{eqnarray*}
and using (\ref{restr}) completes the proof of the lemma

\begin{lemma}\label{lem9} For any $k$ we have
$$|A_k|\geq\frac38\min(|A|^k,(p-1)/2).$$
\end{lemma}
{\bf The proof of Lemma \ref{lem9}.} Let
$$U=\min(|A|^k,(p-1)/2).$$
Then,
$$U-\frac54\frac{U^2}{p-1}\geq\frac38 U,$$
and it suffices to use Lemma \ref{lem8}.

\begin{lemma}\label{lem10} If
$$2\leq k\leq1+\frac{\log((p-1)/2)}{\log|A|}$$
then
$$|N_kA^k|\geq\frac38|A|^{k-8/7}.$$
\end{lemma}
{\bf The proof of Lemma \ref{lem10}.} We use induction on $k$. For $k=2$ the
assertion is trivial. Let us prove it for $k=k+1$
assuming its validity for $k$. The supposition
$k+1\leq1+\log((p-1)/2)/(\log|A|)$ can be rewritten as
$|A|^k\leq(p-1)/2$. By Lemma \ref{lem9}, we have
$$|N_kA^k-N_kA^k|\geq\frac38|A|^k.$$
Applying Corollary \ref{cor7}, we obtain
$$
|4N_kA^k|\geq\left(\frac38|A|^{k-8/7}\right)^{1/8}\left(\frac38|A|^k\right)^{7/8}
=\frac38|A|^{k-1/7}.
$$
Therefore,
$$
|N_{k+1}A^{k+1}|\geq|N_{k+1}A^k|
=|(4N_k+1)A^k|\geq|4N_kA^k|\geq\frac38|A|^{k-1/7},
$$
as required.

\section{The proof of Theorem \ref{th1}} \label{sec6}

We consider several cases.

Case 1. $|A|\leq4$. The inequalities $|A|>p^\delta>p^{1/n}$
show that $4^n>p$ and $|A|\geq2$. Clearly, $|A^n|\geq2$.
By Cauchy--Davenport theorem, for any $N\geq4^n$ we have
$$|NA^n|\geq\min(|A|N-N+1,p)\geq\min(4^n+1,p)=p,$$
and the theorem follows.

Case 2. $4A=\zp$. Then the assertion of the theorem is trivial.

Thus, we assume that (\ref{restr}) holds and $4A\neq\zp$. Using
Cauchy--Davenport theorem, we conclude that
\begin{equation}\label{4A}
|4A|\geq4|A|-3>3|A|.
\end{equation}
Denote
$$n_0=\left[\frac{\log((p-1)/2)}{\log|A|}\right].$$
Notice that $n_0\geq1$ by (\ref{4A}). Also, $n_0\leq n-1$
since
$$|A|^n>p>(p-1)/2\geq|A|^{n_0}.$$

Case 3. Conditions (\ref{restr}), (\ref{4A}) hold and $n_0=n-1$.
By Lemma \ref{lem9} and Lemma \ref{lem10}, for $n\geq 3$ we have
$$|N_{n-1}A^{n-1}-N_{n-1}A^{n-1}|\geq\frac38|A|^{n-1},$$
$$|N_{n-1}A^{n-1}|\geq\frac38|A|^{n-15/7}.$$
These inequalities hold also for $n=2$ if we define $A_1=1$.
Let
$$k=[\log(1/\ve)/\log2]+3.$$
By Corollary \ref{cor7},
$$|kN_{n-1}A^{n-1}|\geq\frac38|A|^{n-1-2^{4-k}/7}
\geq\frac38|A|^{n-1-\ve}.$$
This inequality and (\ref{4A}) imply
$$|kN_{n-1}A^{n-1}||4A|>\frac98|A|^{n-\ve}>p,$$
and we are in position to use Lemma \ref{lem2}:
$$16(kN_{n-1}A^{n-1})(4A)=\zp.$$
Thus,
$$64kN_{n-1}A^n=\zp.$$
We observe that $N_{n-1}\ll 4^n$, $k\ll\log(2+1/\ve)$, and the theorem follows.

Case 4. Conditions (\ref{restr}), (\ref{4A}) hold and $n_0<n-1$.
The last inequality means that $|A|^{n-1}>(p-1)/2$.
By Lemma \ref{lem9} and Lemma \ref{lem10}, we have
$$|N_{n-1}A^{n-1}-N_{n-1}A^{n-1}|\geq\frac3{16}(p-1),$$
$$|N_{n-1}A^{n-1}|\geq|N_{n_0}A^{n_0}|\geq\frac38|A|^{n_0-8/7}
>\frac3{16}(p-1)|A|^{-15/7}.$$
By Corollary \ref{cor7},
$$|3N_{n-1}A^{n-1}|\geq\frac3{16}(p-1)|A|^{-15/{28}}.$$
This inequality and (\ref{4A}) imply
$$|3N_{n-1}A^{n-1}||4A|>\frac9{16}(p-1)|A|^{13/28}.$$
By (\ref{restr}) and (\ref{4A}),
$$|A|^{13/28}>2,\quad p-1>\frac89p.$$
Therefore,
$$|3N_{n-1}A^{n-1}||4A|>p,$$
and we are in position to use Lemma \ref{lem2}:
$$16(3N_{n-1}A^{n-1})(4A)=\zp.$$
Thus,
$$192N_{n-1}A^n=\zp.$$
This completes the proof of the theorem.

\end{document}